%
%

\magnification=1200


\font\AAA=cmr14 at 12pt
\font\BBB=cmr14 at 8pt

\overfullrule=0in

\def\M{{\bf M}}

\def\dim{{\rm dim}}

\def\deg{{\rm deg}}

\def\log{{\rm log}}

\def\arr{\longrightarrow}
\def\supp{{\rm supp}\,}
\def\Link{{\rm Link}}
\def\Wind{{\rm Wind}}


\def\Theorem#1{\medskip\noindent {\AAA T\BBB HEOREM \rm #1.}}
\def\Prop#1{\medskip\noindent {\AAA P\BBB ROPOSITION \rm  #1.}}
\def\Cor#1{\medskip\noindent {\AAA C\BBB OROLLARY \rm #1.}}
\def\Lemma#1{\medskip\noindent {\AAA L\BBB EMMA \rm  #1.}}

\def\Note#1{\medskip\noindent {\AAA N\BBB OTE \rm  #1.}}
\def\Def#1{\medskip\noindent {\AAA D\BBB EFINITION \rm  #1.}}

\def\Conj#1{\medskip\noindent { \AAA C\BBB ONJECTURE \rm    #1.}}

\def\pf{\medskip\noindent {\bf Proof.}\ }
\def\qed{\hfill  $\vrule width5pt height5pt depth0pt$}

\def\w{\wedge}

   \def\cp{{\cal P}}
   \def\co{{\cal O}}
\def\ce{{\cal E}}   
\def\ch{{\cal H}}
\def\cs{{\cal S}}
\def\cd{{\cal D}}

\def\cp{{\cal P}}

\def\and{\qquad {\rm and} \qquad}
\def\arr{\longrightarrow}
\def\ol{\overline}

\def\bbc{{\bf C}}

\def\bbz{{\bf Z}}
\def\bbp{{\bf P}}

\def\a{\alpha}
\def\b{\beta}

\def\e{\epsilon}

\def\g{\gamma}

\def\o{\omega}

\def\s{\sigma}
\def\x{\xi}

\def\L{\Lambda}
\def\G{\Gamma}
\def\O{\Omega}

\def\PH#1{\widehat {#1}}
\def\psh{ \cp\cs\ch_{\o}}

\def\wt{\widetilde}

\def\PL{3}

\def\PHC{5}
\def\PAW{6}
\def\GPAW{2}
\def\PM{3}
\def\MH{4}

\  \vskip .2in
\centerline{\bf   PROJECTIVE LINKING  }
\smallskip

\centerline{\bf
 AND BOUNDARIES OF POSITIVE HOLOMORPHIC CHAINS }
\smallskip

\centerline{\bf IN PROJECTIVE MANIFOLDS, PART II  }
\vskip .2in
\centerline{\bf F. Reese Harvey and H. Blaine Lawson, Jr.$^*$}
\vglue
2cm\smallbreak\footnote{}{ $ {} \sp{ *}{\rm Partially}$  supported by
the N.S.F. }

\centerline{\bf Abstract} \medskip
  \font\abstractfont=cmr10 at 10 pt

{{\parindent=.5in\narrower\abstractfont \noindent 
Part I   introduced the  notion of the {\sl projective linking number }
Link$_\bbp(\G,Z)$ of a compact oriented  real  submanifold $\G$ of dimension  $2p-1$ in complex projective $n$-space $\bbp^n$ with an algebraic subvariety 
$Z\subset \bbp^n-\G$ of codimension $p$. 
It is shown here that a basic conjecture concerning the projective hull of
real curves in $\bbp^2$ implies the following result:   \vskip .1in

\noindent{\sl
$\G$ is the boundary of a positive holomorphic
$p$-chain of mass $\leq \Lambda$ in $\bbp^n$ if and only if the  
$\widetilde{\Link}_\bbp(\G,Z) \geq -\Lambda$ for all algebraic subvarieties 
$Z$ of codimension-$p$ in 
$ \bbp^n-\G$.}

\vskip .1in\noindent
where $\widetilde{\Link}_\bbp(\G,Z) = {\Link}_\bbp(\G,Z)/p!\,\deg(Z)$.
An analogous result is implied in any  projective manifold $X$.

 }}

 \vskip .5in

\centerline{\bf Table of contents.}
\vskip .3in

\moveright .6in\vbox{
1. Introduction.

2. The Projective Alexander-Wermer Theorem.
  
3.  Theorems for General Projective  Manifolds.

4. Variants of the main hypothesis.

}

\vfill\eject

\centerline{\bf 1. Introduction}
\medskip
In Part I we introduced  a    linking  pairing for certain  cycles in projective space
as follows.   Suppose that $\G\subset \bbp^n$ 
  is a compact oriented  submanifold of dimension $2p-1$, and let $Z\subset \bbp^n-\G$ be an algebraic subvariety of codimension $p$.  The {\bf projective linking number} of $\G$ with $Z$ is defined to be
 $$
  \Link_{\bbp}(\G, Z) \  \equiv \ N\bullet Z - \deg (Z) \int_N\omega ^p
  $$
where $\omega$ is the standard  K\"ahler form on $\bbp^n$ and $N$ is any  integral 2p-chain
with $\partial N=\G$ in $\bbp^n$. This definition is independent of the choice of $N$. 
  The associated {\bf reduced linking number} is defined to be
  $$
  \widetilde{\Link}_{\bbp}(\G,Z) \ \equiv \ {1\over p!\, \deg(Z)}\Link_{\bbp}(\G,Z)
  $$
  The basic result proved here is the following.
  
  \Theorem{1.1} {\sl  Let $\G\subset \bbp^n$ be a compact oriented real analytic submanifold of dimension $2p-1$ with possible integer multiplicities on each component. If  Conjecture B holds
   (See below), then the following are equivalent:
  \medskip
  
  (1)\ \ \ $\G$ is the boundary of a positive holomorphic $p$-chain of mass $\leq  \Lambda$ in  $\bbp^n$.
  \medskip
  
(2)\ \  $ \widetilde{\Link}_{\bbp}(\G,Z)\ \geq \   -\Lambda$ for all algebraic subvarieties $Z$ of         codimension $p$ in $\bbp^n-\G$.
}
\medskip  
  
 A compact subset $K\subset \bbp^n$ is called {\bf stable} if the best constant function is bounded on the projective hull $\PH K$ (See [HL$_{3,4}$]).
It is known that  for any stable real analytic curve $\gamma \subset \bbp^n $, the set 
  $\PH \g-\g$ is a 1-dimensional complex analytic subvariety of $\bbp^n-\g$, [HLW].
  Conjecture A from Part I is the statement that any compact real analytic curve  $\gamma \subset \bbp^n $ is stable.  Even more likely is the following.
  
  \Conj{B} Let $\gamma \subset \bbp^2$ be a compact embedded real analytic curve such that for some
  choice of orientation and positive integer multiplicity on each component, condition (2) above is
  satisfied.  Then $\g$ is stable.

\medskip

In Part I the conclusion of Theorem 1.1 was established for any stable real analytic curve $\G\subset
\bbp^n$ (with orientation and multiplicity on each component). The main point of   Part II is to
prove that this result for $p=1$ implies the result for all $p>1$, provided one can drop the stability hypothesis in the $p=1$ case.   

Note incidentally that there is no assumption of maximal complexity on the cycle $\G$.

Theorem 1.1 represents a projective analogue of a result of H. Alexander and J. Wermer [AW].

If the cycle $\G$ in Theorem 1.1 bounds a holomorphic $p$-chain $T$, then there is a unique such
chain $T_0$ of least mass with $dT_0=\G$. (All others are obtained by adding positive algebraic 
$p$-cycles to $T_0$.)  

\Cor{1.2}  {\sl Let $\G$ be as in Theorem 1.1 and suppose that $\G=dT$ for some positive holomorphic $p$-chain $T$.  Then $T$ is the chain of least mass with boundary $\G$ if and only if
$$
\inf_Z \left\{ {  T\bullet Z \over \deg Z  }  \right\} \ =\ 0
$$
where the infimum  is taken over all positive algebraic $(n-p)$-cycles in $\bbp^n-\G$.}

\medskip

The linking hypothesis (2) in Theorem 1.1 can be replaced by other hypotheses.  This is discussed
in \S 3. Another interesting consequence of Theorem 1.1 is the following result, whose proof 
follows exactly the lines given in Part I for the case $p=1$.

\Theorem {1.3}  {\sl  Let $M\subset \bbp^n$ be a compact  embedded real analytic submanifold of dimension $2p-1$ and assume Conjecture B.  Then a class $\tau \in H_{2p}(\bbp^n, M; \bbz)$
contains a positive holomorphic chain $T$ (with $\supp dT \subseteq M$) if and only if 
$$
\tau\bullet u\ \geq\ 0
$$
for all classes $u\in H_{2n-2p}(\bbp^n - M;\bbz)$ which are represented by positive algebraic cycles.}
\medskip

Theorem 1.1 and many of its consequences carry over to general projective manifolds.  
This is done in \S 4.
  
We recall our convention that  
$
d^C\ =\ {i\over 2\pi}(\overline\partial - \partial)
$


\vskip .3in

\centerline{\bf 2. The Projective Alexander-Wermer Theorem.} 
\medskip

Let $\G$ be a compact smooth oriented submanifold of dimension $2p-1$  in $\bbp^n$.  We recall that
(even if $\G$ is only class $C^1$) any irreducible complex analytic
subvariety $V\subset \bbp^n-\G$ has finite Hausdorff 2$p$-measure and defines
a current $[V]$ of dimension 2$p$ in $\bbp^n$ by integration on the
canonically oriented manifold of regular points.  Furthermore, the boundary of this current is of the form $d[V] =\sum_j \epsilon_j \G_j$ where $\G_1,...,\G_\ell$
represent the connected components of $\G$ and $\epsilon_j =1, 0, $ or $-1$.
(See [H] for example.)  We now allow $\G$ to carry positive integer multiplicities on each 
component, so it is of the form $\G =\sum_j m_j \G_j$.

\Def{\GPAW.1} By a {\bf positive holomorphic $p$-chain with boundary $\G$} we
mean a finite sum $T=\sum_k n_k [V_k]$ where each $n_k\in\bbz^+$ and
each $V_k\subset \bbp^n-\G$ is an irreducible subvariety of dimension $p$, so
that
$$
dT\ =\ \G\qquad\ \ {\rm (as\ currents\ on }\ \  \bbp^n)
$$

\medskip

By the {\bf mass} of such a chain $T=\sum_k n_k [V_k]$ we mean its weighted volume:
$\M(T) \equiv \sum_k n_k \ch^{2p}(V_k)= T(\O_p)$ where $\ch^{2p}$ denotes Hausdorff
2$p$-measure and
$$
\O_p \ \equiv \ {{1\over p!}}\o^p
$$

\Prop{\GPAW.2}  {\sl Suppose $T$ is a positive holomorphic $p$-chain with boundary $\G$ as above.
Then 
$$
 \widetilde{\Link}_{\bbp}(\G,Z) \ \geq \ -\M(T)
$$
for all positive algebraic cycles $Z$ with support in $\bbp^n-\G$.}

\pf
Note that  since $dT=\G$ we have 
$$
 \widetilde\Link_{\bbp}(\G, Z) = {T\bullet Z \over p!\,\deg Z}-  T\left(\O_p\right) \geq 
 -T\left(\O_p\right)=  -\M(T)
 $$
  since $T\bullet Z\geq 0$ by the positivity of $T$ and $Z$.
\qed

\medskip
Note that Proposition 2.2 holds for positive holomorphic chains with quite general boundaries $\G$.
This brings us to the main result.

\Theorem {\GPAW.3} {\sl Under the assumption of  Conjecture B
the following are equivalent.
\medskip

(1) \ $\G=dT$
where $T$ is a positive holomorphic $p$-chain in $\bbp^n$ with $\M(T)\leq
\Lambda$
\medskip

(2)\ \ 
$
\widetilde{\Link}_{\bbp}(\G,Z)\ \geq\ -\Lambda
$
for all $(n-p)$-dimensional algebraic varieties $Z\subset \bbp^n-\G$.}

\pf  Proposition \GPAW.2 states that  (1) $\Rightarrow$ (2).
For the converse we shall show that the linking condition
persists for hyperplane slices, and then proceed by induction on dimension.

\Prop{\GPAW.4} {\sl
Suppose that $\G$ satisfies the $\Lambda$-linking condition (2) in $\bbp^n$.
If $H \cong\bbp^{n-1}$ is a hyperplane which intersects $\G$ transversely,
then $\G_H\equiv \G\cap H$ satisfies the $\Lambda'$-linking condition in $H$
where 
$$
\Lambda'\ =\  p\Lambda + \int_\G d^C u\wedge\O_{p-1}
$$
and $u=\log(|Z_0|/\|Z\|)$ where $Z_0$ is the linear function
defining $H$. }

\pf Since bordism and homology agree in $\bbp^n$ there exists a
 compact oriented $2p$-manifold $N$ with 
boundary and a smooth map $f:N\to \bbp^n$ such that
$f$ is an immersion near $\partial N$ and 
$$
f\bigl|_{\partial N}:\partial N \to \G \qquad {\rm is\ an \ oriented\
diffeomorphism}.
$$
 Since $\G$ is transversal to $H$, $f$ is 
also transversal to $H$ near the boundary.  By standard transversality
theory we can  perturb $f$, keeping it fixed near the boundary,           
so that it is everywhere transversal to $H$. Let $N_H \equiv f^{-1}(H)$
oriented by $N$ and the normal bundle to $H$, and let $f_H:N_H\to H$ be the
restriction of $f$.  Then $(f_H)_*[N_H]$ defines a $(2p-2)$-dimensional
current in $H$ with boundary $\G_H$.  We denote this current simply by 
$[N_H]$.

Suppose now that $Z\subset H-\G_H$ is an $(n-p)$-dimensional algebraic
subvariety.  We may assume, again by a small perturbation, that
$f$ misses the singular set of $Z$ and is transversal to Reg$(Z)$. 
It is then straightforward to check that 
$$
[N]\bullet Z\ =\ [N_H]\bullet_H Z
\eqno{(\GPAW.1)}
$$
where ``$\bullet_H$'' denotes the intersection pairing in $H$ (defined as
in \S\PL \ of Part I).

By assumption we have that
$$
\widetilde{\Link}_{\bbp}(\G,Z) \ =\   {1\over p!}  \left\{  { 1\over{\deg Z}} ([N]\bullet Z)
-\int_N\o^p  \right\} \ \geq\ -\Lambda.
\eqno{(\GPAW.2)}
$$
Now the function $u$ above satisfies the Poincar\'e-Lelong equation
$$
dd^C u\ =\  H - \omega.
\eqno{(\GPAW.3)}
$$
Substituting (\GPAW.1) and (\GPAW.3)  into (\GPAW.2) gives
$$\eqalign{
{ 1\over{\deg Z}} ([N_H]\bullet_H Z) - &\int_N (H-dd^Cu )\w \o^{p-1}
\ = \cr
&\qquad = \ { 1\over{\deg Z}} ([N_H]\bullet_H Z) - \int_{N_H}\o^{p-1} 
+  \int_N dd^C u\w \o^{p-1}   \cr
&\qquad = \ \widetilde{\Link}_{\bbp}(\G_H,Z) (p-1)! + \int_\G d^C u \w \o^{p-1}  \cr
&\qquad \geq \ -\Lambda p!
}
$$ 
where the first equality is straightforwardly justified using
transversality.\qed

\Cor{\GPAW.5} {\sl  Assume that (2) $\Rightarrow$ (1) for all manifolds
$\G$ of dimension $2p-3$ in projective space.  Suppose that $\G$ is a
$(2p-1)$-manifold satisfying (2) and that $H$ is a hyperplane transversal
to $\G$. Then there exists $\Lambda'>0$ so that $\G_{H'}\equiv \G\cap H'$ bounds
a positive holomorphic $(p-1)$-chain of mass $\leq \Lambda'$ for all
hyperplanes $H'$ in a neighborhood $U$ of $H$. }

\pf
If $H$ is transversal to $\G$ then so are all hyperplanes $H'$ in a
neighborhood of $H$. Furthermore, the integral 
$
\int_\G d^C (u_{H'})\w\o^{p-1}
$
depends continuously on $H'$ in that neighborhood, where
$u_{H'} =  \log(|(A_{H'}, Z)|/\|Z\|)$
and $A_{H'}$ is a continuous choice of vectors with
$H'=\{[Z]\in\bbp^n: (A_{H'},Z)=0\}$. It follows that the constant
$\Lambda'$ in Proposition \GPAW.4 is uniformly bounded below in a neighbohhood of
$H$.  One then applies the inductive hypothesis.\qed

\medskip

\Prop{\GPAW.6}  {\sl Let $U$ be the neighborhood given in Corollary
{\GPAW.5}. For each hyperplane $H'$ in $U$ let $T_{H'}$ be the
positive holomorphic $p-1$ chain {\sl of least mass} with $dT_{H'} = \G_{H'}$.
Then $T_{H'}$ is uniquely determined by $H'$ and the mapping 
$H'\mapsto T_{H'}$ is continuous on $U$.}

\pf
We first prove uniqueness.  For future reference we formulate this result for $\G$ instead of $\G_{H'}$.

\Lemma{\GPAW.7} {\sl  Suppose $\G$ bounds a positive holomorphic $p$-chain.  Then 
the positive holomorphic $p$-chain of least mass with boundary $\G$ is unique.}

\pf
Suppose that $T=\sum_i n_i [V_i]$ and 
$T'=\sum_j n'_j [V'_j]$ are
positive holomorphic $p$-chains of least mass having the same boundary 
$dT=dT'=\G$ in $\bbp^n$. By the least mass hypothesis we know that 
$$
d[V_i]\ \neq\ 0\and d[V'_{j}]\ \neq\ 0\qquad {\rm  for\ all\ \ } i,j.
\eqno{(\GPAW.4)}
$$
In fact $d[V_i]$ and $d[V'_j]$ each consist of a finite number of 
oriented connected components of $\G$, each  with multiplicity one.
Let $\G_1$ be an oriented connected component of $\G$. (If $\G$ has  
 multiplicity greater than 1 along $\G_1$, we ignore that multiplicity for
the moment.)  Then there must exists a component of $T$, say $V_1$, such
that  $\G_1$ forms part of the oriented boundary $dV_1$. Similarly there is a
component, say $V'_1$ of $T'$ such that $\G_1$ is part of $dV'_1$.
By boundary regularity [HL$_1$], and local and global uniqueness these two 
irreducible subvarieties of $\bbp^n-\G$ must coincide.  Hence, $S
\equiv T-V_1$   and $S'\equiv T'-V'_1$ are positive holomorphic $p$ chains
with $dS=dS'$. By continuing this process one of the two chains will
eventually be reduced to zero. However, the other must also be 0 since its
boundary is zero and its remaining components satisfy condition (\GPAW.4)
\qed\medskip

To prove continuity  it will suffice to show that every convergent sequence
 $H_j \to H$  has a subsequence  such that $T_{H_j}\to T_H$. By the local
uniform bound on the mass, the fact that $dT_{H_j}\to dT_H=\G_H$, and the
compactness of positive holomorphic chains, we know that there is a
subsequence which converges to some  positive holomorphic chain $T$ with
boundary $\G_H$. We then apply the uniquness.\qed

\medskip

We recall that $\G$ is said to be {\bf maximally complex} if 
$$
\dim_{\bbc}(T_x\G \cap JT_x\G) \ =\ p-1  \qquad
{\rm for\ all\ \ }  x\in \G
$$
where $J$ is the almost complex structure on $\bbp^n$.

\Prop{\GPAW.8}  {\sl Assume Conjecture B. If $\G$ satisfies (2), then $\G$ is
maximally complex. }

\pf  
The result is trivial when $\dim \g =1$ so we first consider the case
$\G= \G^3\subset \bbp^3$.  We want to show that   $\int_\G \a=0$ for all 
$(3,0)$-forms $\a$ on $\bbp^3$.

Choose a line $L\cong \bbp^1\subset\bbp^3$ with
$\G\cap L = \emptyset$ and a linear projection $\pi:\bbp^3-L\to \bbp^1$.
Fix a point $x_\infty\in \bbp^1$ and choose affine coordinates
$(z_0,z_1,z_2)$ on $\pi^{-1}(\bbp^1-\{x_\infty\})\cong\bbc^3$.
We shall show that 
$$
\int_\G g(z_0) dz_0\w dz_1\w dz_2 =0
\eqno(\GPAW.5)
$$
 for all $g\in C_0^\infty(\bbc)$. Such forms, taken over all possible
choices above, are dense  in $\ce^{3,0}$ on a neighborhood of $\G$.
Hence,  $\int_\G \a=0$ for all  $(3,0)$-forms $\a$, which implies
that $\G$ is maximally complex.

To prove (\GPAW.5) we choose a 4-chain $N$ with compact support in 
$\bbp^3-L$ such that $d N=\G$. Then (\GPAW.5) can be rewritten as
$$
\int_N{{\partial g}\over{\partial \ol z_0}}\ d z_0 \w d{\ol z}_0\w dz_1\w
dz_2 \ =\ 0
$$ 
for all $g\in C_0^\infty(\bbc)$. For this it will suffice to prove that
$$
(N\w dz_1\w dz_2, \pi^*  \eta)\ =\ (\pi_* (N\w dz_1\w dz_2),  \eta)\ =\  0
$$
for any  (1,1)-form $\eta$ with compact support in $\bbc$. For this it will
suffice to consider $\eta = \delta(z_0-t) dz_0\wedge d{\ol z}_0$ for $t\in
\bbc$, in other words we want to show that the slice at $t$:
$$
\{\pi_*(N \w dz_1\w dz_2)\}_t \ =\ \pi_*\{N_t\w  dz_1\w dz_2\}\ =\
\int_{N_t} dz_1\w dz_2\ =\ 0 
\eqno(\GPAW.6)$$
for all $t\in \bbc$.  

Observe now that $d N_t \ =\ \G_t$, the slice of $\G$ by $\pi$
at $t$, and by Proposition \GPAW.4 this $\G_t$ satisfies the projective
linking condition (2). Hence by Theorem \PAW.1 in [HL$_4$] and our hypothesis
that $\G_t$ is stable, we conclude that $\G_t= dT_t$
where $T_t$ is a positive holomorphic 1-chain in $\bbp^2$ = the closure of
$\pi^{-1}(t)$. Thus our desired condition (\GPAW.6) is established by the
following.

\Lemma{\GPAW.9} {\sl Let $\g$ be a curve, or in fact any rectifiable
1-cycle with compact support in $\bbc^2\subset \bbp^2$. Suppose $\g = dT$
where  $T$ is a  positive holomorphic chain in $\bbp^2$.  Then for any  
$S\in\cd'_{2, {\rm cpt}}(\bbc^2)$ with $dS=\g$, one has $S(dz_1\w dz_2)=0$.}

\pf It suffices to construct one current $S$ with these properties.
We can assume that the line at infinity $\bbp^1_{\infty}=\bbp^2-\bbc^2$
meets $\supp T$ only at regular points and is transversal there.
(The general result follows directly.) Choose $x\in \supp\, T \cap
\bbp^1_\infty$ and let $L$ be the tangent line to $\supp\, T$ at $x$.
Then after an affine transformation of the $(z_1,z_2)$-coordinates,
 we may assume $L\cong z_1$-axis. This transformation can be chosen
with determinant one, so the form $dz_1\w dz_2$ remains unchanged.
 
Near the point $x$, the current $T$ is given by a positive  multiple of 
the graph 
$
\Sigma_R\ \equiv\  \{(z_1, f({1 \over{z_1}})) : |z_1|\geq R\}
$
where $f$ is holomorphic in the disk of radius $1/R$ and satisfies
$$
f(0)\ =\ f'(0)\ =\ 0.
$$
In particular we have that
$$
\lim_{z_1\to 0}z_1f(1/z_1)\ =\ 0.
\eqno{(\GPAW.7)}$$
We now modify $T$ by replacing (the appropriate multiple of) $\Sigma_R$
with the current $L_R + U_R$ where 
$$
L_R\ \equiv\ \{(z_1,0) : |z_1|\leq R\} \and
U_R\ \equiv\ \{(z_1,tf(1/z_1) ) : |z_1|= R \ \ {\rm and\ \ } 0\leq t\leq 1\}
$$
with orientations chosen so that $d(L_R + U_R)=d\Sigma_R$.
Observe that 
$$
\int_{L_R + U_R} dz_1\w dz_2 \ =\ \int_{U_R} dz_1\w dz_2 \ =\ 
-\int_{dU_R} z_2 dz_1\ =\
\int_0^{2\pi}f(e^{i\theta}/R){{ie^{i\theta}}\over{R}} d\theta\  \to\  0
 $$
as $R\to \infty$ by (\GPAW.7).

Carrying out this procedure at each point of $\supp T\cap \bbp^1_{\infty}$
we obtain a current $T(R_1,...,R_\ell)$ with compact support in $\bbc^2$
and with $dT(R_1,...,R_\ell)=\g$.  Since $dz_1\w dz_2$ is closed we have
$$
S(dz_1\w dz_2)\ =\ T(R_1,...,R_\ell) (dz_1\w dz_2)\ =\ \sum_{k=1}^\ell
\int_{ U_{R_k}} dz_1\w dz_2\ \to\ 0
$$
as $R_1,...,R_\ell\to \infty$.
\qed\medskip

We have now proved the proposition for 3-folds in $\bbp^3$. The result for 
3-folds in $\bbp^n$ follows by considering the family of projections
$\bbp^n - - - > \bbp^3$ which are well defined on $\G$. The result for
general $\G$ follows by intersecting with hyperplanes and applying 
Proposition \GPAW.4.  \qed  \medskip

We now show that $\G$ bounds a positive holomorphic chain by applying the main result
in [DH]. Let $L\subset \bbp^n$ be a linear subspace of
(complex) codimension $p-1$ which is transversal to $\G$.
We  assume that $L$ meets every component of $\G$. This can
be arranged by taking a Veronese embedding   $\bbp^n\subset \bbp^N$ 
of sufficiently high degree, and working with linear subspaces there.
One checks that the projective linking numbers of $\G$ are also bounded 
in $\bbp^N$.  By applying Proposition \GPAW.4 inductively we see that
for all linear subspaces $L'$ in a neighborhood of $L$, the intersections
${\G}_{L'} = \G\bullet L'$ satisfy the projective linking condition for oriented 
curves with multiplicities.  Therefore by  [HL$_4$] and our assumption of 
Conjecture B, each slice  ${\G}_{L'}$ bounds a positive holomorphic 1-chain.
With this property and maximal complexity, 
It follows directly from [DH] that $\G$ bounds a holomorphic $p$-chain.
The unique minimal such chain $T$ will be supported in the subvariety
$\bbp^n\subset \bbp^N$ because $\G$ is. Furthermore, $T$ must be 
positive.  If not, there would be a negative component, say $T_0$. 
Since $T$ is minimal, $dT_0\neq 0$. Now $L$ must meet $T_0$ since
it meets all components of $\G$. It follows that the minimal holomorphic
1-chain with boundary ${\G}_L$ is not positive -- a contradiction. 

\Note{}This last paragraph could be replaced by an argument based on the
 results in [HL$_2$].
\medskip
 
 So we have proved that $\G=dT$ where $T$ is a positive holomorphic $p$-chain.
We may assume that $T$ is the unique such chain of least mass.  It remains to prove
that $\M(T)\leq \Lambda$. 

Suppose not.  Then 
$$
\M(T) \ =\ T(\O_p)\ = \ \Lambda +r
\eqno{(\GPAW.8)}
$$
for  $r>0$, and we see that 
$$
\eqalign{
{ T\bullet Z \over p!\,\deg Z } \ &= \ { T\bullet Z \over p!\,\deg Z } -T(\O_p)+T(\O_p) \cr
&=\   \widetilde{\Link}_{\bbp}(\G,Z) + T(\O_p)  \ \geq\  - \Lambda + \Lambda +r=r\cr
}
$$
for all algebraic subvarieties $Z$ of codimension $p$ in  $\bbp^n-\G$.
Hence, it will suffice to prove that 
$$
\inf_Z\left\{  { T\bullet Z \over \,\deg Z }  \right\}\ =\ 0
\eqno{(\GPAW.9)}
$$
where $Z$ varies over the codimension $p$ subvarieties of $\bbp^n-\G$.

\Lemma{\GPAW.10}  {\sl Let $T$ be a positive holomorphic $p$-chain in $\bbp^n$ with
boundary $\G$.  Then $T$ is the unique such chain of least mass if and only if 
every irreducible component of $\supp T$ in $\bbp^n-\G$ has a non-empty boundary
(consisting of components of $\G$).}

\pf  If $T$ has least mass, it can  have no components with boundary zero, since one could
remove these components and thereby reduce the mass without changing the boundary.
If $T$ is not of least mass, let $T_0$ be the least mass solution and note that
$d(T-T_0)=0$. It follows from [HS] (actually, elementary arguments involving local uniqueness at the boundary will suffice here) that $T-T_0 = S-S_0$ where $S$ and $S_0$ are positive algebraic
$p$-cycles, i.e., $dS=dS_0=0$. We assume that $S$ and $S_0$ have no components in common.
From the equation $T+S_0=T_0+S$ and the uniqueness of the decomposition of analytic subvarieties
into irreducible components, we see that the components of $S$  must be components of $T$ and 
similarly the components of $S_0$  must be components of $T_0$. Since $T_0$ is least mass, we have $S_0=0$ and $T=T_0+S$.  Since $T$ is not the least mass solution, $S\neq 0$.\qed

Our proof now proceeds by induction on $p$.  We assume that condition (\GPAW.9) holds for least mass
chains of  dimension $<p$, and we shall show that it hold in dimension $p$.  The case
$p=1$ has already been established in Part I, Corollary 6.8.

  We return to our  positive holomorphic $p$-chain $T$ of least mass. A positive divisor $D$
  in $\bbp^n$ is defined to be {\bf totally transverse} to $T$ if: 
  
  (i)  $D$ is smooth, 
  
  (ii) $D$ meets every component of $dT=\G$, and each of these intersections is transverse,
  
  (iii)  $D$ is transversal to every stratum of the singular stratification of $\supp T$.
  
  \noindent
  This condition is open,  and it is non-empty for divisors of sufficiently high degree.
  Let $D$ be such a divisor with degree $d$.  Let $\bbp^n\subset \bbp^N$ be the 
  order $d$ Veronese embedding, and let $H_0\subset \bbp^N$ be the hyperplane with
  $D=H_0\bullet \bbp^n$.  Then $H_0$ is totally transverse to $T$ in $\bbp^N$.
  Assume there exists $H$ in a neighborhood of $H_0$ which is totally transverse to $T$ and such that
  $T_{H}\equiv H\bullet T$ has no irreducible components with boundary zero.  Then there exists
  a sequence $\{Z_j\}_{j\geq 0}$ of subvarieties of codimension $p-1$ in $H-\G$ such that
  $$
  \lim_{j\to\infty} \left\{  { T_H\bullet_H Z_j \over \,\deg Z_j }  \right\}\ =\ 0.
  $$
  Note that $T_H\bullet_H Z_j = T\bullet_{\bbp^N} Z_j$, and so
   $$
  \lim_{j\to\infty} \left\{  { T\bullet_{\bbp^N}  Z_j \over \,\deg Z_j }  \right\}\ =\ 0.
  $$  
  Now by a small perturbation we may assume that each $Z_j$ is transversal to $\bbp^n$.
  Then $W_j\equiv Z_j\cap \bbp^n$ is a subvariety of codimension $p$ in $\bbp^n$ with 
   $T\bullet_{\bbp^n}  W_j =T\bullet_{\bbp^N}  Z_j $ and with   degree $\deg W_j = d^p \deg Z_j$.
  It follows that 
   $$
  \lim_{j\to\infty} \left\{  { T\bullet_{\bbp^n}  W_j \over \,\deg W_j }  \right\}\ =\ 0
  $$  
as desired. 

 So we are done unless every  totally transverse hyperplane section  $T_H$   has components with
 zero boundary. We will show that this cannot happen for $H$ in a neighborhood of $H_0$.
 For the sake of clarity, we assume first that each component of $\supp T$ is smooth. Then by transversality each component of $\supp T_H$ is smooth.  Suppose there exists a component $V$  of $\supp T_H$ which is without boundary. Then $V$ will be contained in one of the 
 components of $\supp T$, which we denote by $W$. Note that $W$ is a smooth submanifolds
 with real analytic boundary.
 
 Now let $V_\epsilon$ be a neighborhood of 
 $V$ in  $W$. Since $V_\epsilon$ is a subvariety defined in a neighborhood of a 
 hyperplane, it extends to an irreducible {\bf algebraic} subvariety $Y$ of dimension $p$.
It follows that $W$ is a subdomain with real analytic boundary in $Y$.

Now the generic hyperplane section of an irreducible variety is again irreducible when
$p>1$ (See [Ha, Prop. 18.10]).  It follows that $W \cap H = V$.  However, this is impossible since $H$ meets every  component of $\G$ and so it must meet the components of $\G$ 
which are contained in $W$.

When $\supp T$ is not smooth the argument is similar. By total transversality  and our assumption,
there exists a component $V$ of $\supp T_H$ (for generic $H$) with no boundary and with the
property that $V$ extends to an irreducible $p$-dimensional subvariety $V_\epsilon$ in a 
neighborhood of $H$. By [HL$_1$, Thm. 9.2],   $V_\epsilon$ extends to an irreducible algebraic subvariety  $Y$ of dimension $p$ in $\bbp^N$. The remainder of the argument is the same.
\qed

\medskip

There is an ``affine'' version of Theorem \GPAW.3 parallel to the 
``affine'' version  (Theorem 6.6 ) of Theorem 6.1 in Part I.  However the reader should 
note that the hypothesis that $\G$ is contained in some affine chart 
is not satisfied generically when $\dim(\G) >1$.

\Theorem{\GPAW.10}  {\sl Let $\G$ is as in Theorem \GPAW.3  and suppose $\G\subset \bbc^n$ for some affine chart  $\bbc^n\subset\bbp^n$. Assuming Conjecture B, the following are equivalent:
\medskip

(1)\ \ There exists a constant $\Lambda$ so that the classical linking number satisfies
$$
{\Link}_{\bbc^n}(\G,Z)\ \geq\ -\Lambda\,p! \,\deg Z
$$

\qquad for all $(n-p)$-dimensional algebraic varieties $Z\subset \bbc^n-\G$.
\medskip

(2) \ \ \ $\G=dT$
where $T$ is a positive holomorphic $p$-chain in $\bbp^n$ with }
$$
\,\M(T)\leq \Lambda + {1\over p}\int_\G d^C\log\sqrt{1+\|z\|^2}\wedge\O_{p-1}.
$$

\pf  We recall that in terms of the affine coordinate $z$ on $\bbc^n$, the Kaehler form is 
given by $\o={1\over 2}dd^c\log(1+\|z\|^2)$.  Now let $N$ be a 2$p$-chain in $\bbc^n$ with $dN=\G$. Then
$ \widetilde\Link_{\bbp}(\G, Z) p!  = {1\over{\deg Z}}N\bullet Z - \int_N \omega^p 
=  {1\over{\deg Z}}\Link_{\bbc^n}(\G, Z) - \int_\G d^C
\log\sqrt{1+\|z\|^2}\wedge \omega^{p-1}$.
The result now follows directly from Theorem \GPAW.3.
\qed


\vskip .3in

\centerline{\bf \PM. Theorems for General Projective Manifolds.}  \medskip
  
The results established above generalize from $\bbp^n$ to any projective
manifold.  Let $X$ be a compact complex $n$-manifold with a positive 
holomorphic line bundle $\lambda$.
Fix a hermitian metric on $\lambda$ with curvature form $\omega>0$, and
give $X$ the K\"ahler metric associated to $\omega$.
Let $\G$ be a  $(2p-1)$-cycle on $X$ with properties as in \S \GPAW\
(i.e., an oriented $(2p-1)$-dimensional submanifold with integral weights),
and assume $[\G]=0$ in $H_{2p-1}(X;\,\bbz)$.

\Def{\PM.1}  Let $Z$ be a positive algebraic $(n-p)$-cycle on $X$ 
which has cohomology class $\ell [\o^p]$ for some $\ell\geq 1$. If $Z$
does not meet $\G$, we can define the  {\bf linking number} and the {\bf
reduced linking number} by
$$
\Link_{\lambda}(\G,Z)\ \equiv \ N\bullet Z-\ell\int_N\omega^p
\and
\widetilde{\Link}_{\lambda}(\G,Z)\ \equiv \ {1\over \ell \, p!}\ \Link(\G,Z)
$$
respectively, where $N$ is any 2$p$-chain in $X$ with $dN=\G$ and where
the intersection pairing $\bullet$ is defined as in  \S 3 of  Part I with $\bbp^n$
replaced by $X$.

\medskip

To see that this is well-defined suppose that $N'$ is another 2$p$-chain
with $dN'=\G$. Then
 $(N-N')\bullet Z - \ell\int_{N-N'}\omega^p = (N-N')\bullet
(Z-\ell[\omega^p])=0$ because $Z-\ell\omega^p$ is cohomologous to zero in $X$.

\Theorem{\PM.2} {\sl  Under the assumption of Conjecture B the following are
equivalent:

\medskip

(1) \ \ $\G=dT$\ \ where $T$ is a positive holomorphic $p$-chain on $X$ with
mass $\leq  \Lambda$.

\medskip

(2) \ \ $\widetilde{\Link}_{\lambda}(\G,Z) \geq -\Lambda$ \ \ for all
positive algebraic $(n-p)$-cycles $Z\subset X-\G$ of

\ \ \ \ \ \ \ cohomology class  $\ell[\omega]^p$ for $\ell\in \bbz^+$.
}

\pf That (1) $\Rightarrow$ (2) follows as in the proof of Proposition 2.2. In fact this
shows that (2) holds with no restriction on the cohomology class of $Z$.

For the converse we may assume (by replacing $\lambda$ with $\lambda^m$
if necessary) that the full space of sections $H^0(X,\co{}(\lambda))$ gives
an embedding $X\subset \bbp^N$. An algebraic subvariety $\wt Z \subset
\bbp^N$ of codimension $p$ is said to be {\sl transversal} to $X$ if each
level of the singular stratification of $\wt Z$ is transversal to $X$.
More generally, a positive algebraic cycle $\wt T=\sum_\a n_\a {\wt Z}_\a$
of codimension $p$ is {\sl transversal} to $X$ if each $Z_\a$ is. Such
cycles are dense in the Chow variety of all positive algebraic
$(N-p)$-cycles in $\bbp^N$. This follows from the Transversality  Theorem for Families   applied to the family GL$_{\bbc}(N+1)\cdot\wt{Z}$ and the submanifold $X$ in $\bbp^N$ (cf. [HL$_1$, App. A]). 

It is straightforward to check that if $\wt Z$ is
transversal to $X$, then $Z = \wt Z \cap X$ is an algebraic
subvariety  of codimension $p$ in $X$ with cohomology class $\ell [\omega^p]$
where $\ell = \deg\ \wt Z$.  Let $\wt Z$ be such a cycle with the property
that $\wt Z \cap \G=\emptyset$. Let $N$ be a 2p-chain in $X$ with boundary
$\G$ which meets $Z=\wt Z \cap X$ transversely at regular points. Then
local computation of intersection numbers shows that $N\bullet_X Z
= N\bullet_{\bbp^N} \wt Z$.   Consequently,   hypothesis (2) implies that 
$$
\widetilde{\Link}_{\bbp}(\G,\wt Z)\ =\ 
\widetilde{\Link}_{\lambda}(\G,Z) \geq -\Lambda.
$$
Since this holds for a dense set of subvarieties of $\bbp^n-\G$ it holds for
all such subvarieties.  Theorem \GPAW.3 now implies that $\G$ bounds a
holomorphic chain $T$ in $\bbp^N$.  Since $\G$ is supported in $X$, so
also is $T$.\qed


\vskip .3in

\centerline{\bf \MH.  Variants of the Main Hypothesis}
\medskip

The linking hypothesis (2) in Theorem \PM.2 can be replaced by several
quite different conditions thereby yielding several geometrically
distinct results. In this section we shall examine these conditions.

Let $X, \lambda, \omega$ and $\G$ all be as in \S \PM. Suppose 
$Z\subset X$ is an algebraic subvariety of codimension-$p$ with cohomology
class  $\ell[\omega^p]$.  Then a {\bf spark} associated to $Z$  with curvature
$\ell\o^p$  is a current $\a\in \cd'^{2p-1}(X)$ which satisfies the
{\sl spark equation}  
$$
d \a\ =\  Z\ -\  \ell\o^p.
\eqno{(\MH.1)}
$$
Such sparks form the basis of (one formulation of) the theory of 
differential characters (cf. [HLZ]). Two sparks $\a, \a'$ satisfying 
(\MH.1) will be called {\bf  commensurate} if $\a'=\a+d\b$ for
$\b \in\cd'^{2p-2}(X)$.

\medskip
\noindent 
{\bf A. $\lambda$-Winding Numbers.}
Suppose now that $\G\cap Z=\emptyset$.  Up to commensurability we
may assume  $\a$ is smooth in a neighborhood of $\G$ (See [HLZ, Prop. 4.2]).
The {\bf  reduced $\lambda$-winding number} can then be defined as
$$
\widetilde{\Wind}_{\lambda}(\G, \a)\ \equiv\   {{1}\over{\ell\, p!}}\int_\G \a.
$$
It follows from the spark equation (\MH.1) that
$$
\widetilde{\Wind}_{\lambda}(\G, \a)\ =\ \widetilde{\Link}_{\lambda}(\G, Z).
\eqno{(\MH.2)}
$$
In particular, this winding number is independent of the choice of $\a$.  
This can be directly verified using deRham theory and the fact that
$[\G]=0$ in $H_{2p-1}(X;\bbz)$.

\medskip
\noindent 
{\bf B. Positivity.}
Recall that a smooth $(p,p)$-form $\psi$ on $X$ is
called  {\bf weakly positive} if 
$
\psi(\x)\ \geq\ 0
$
for all simple $2p$-vectors $\x$ representing canonically oriented complex
tangent $p$-planes to $X$.  A current $T\in \cd'_{p,p}(X)$ of bidimension
p,p is called {\bf  positive} if $T(\psi)\geq0$ for all weakly positive 
$(p,p)$-forms $\psi$  on $X$. In this case we write $T\geq 0$.

\medskip
\noindent 
{\bf C. Algebraic Homology.} Set
 $H^+_{2k} \equiv \{z\in H_{2k}(X-|\G|;\bbz): (\omega^k,z)\geq0\}$, where $|\G|\equiv \supp \G$,  and
define $$H^+_{2k, {\rm alg}} \ \subseteq\  H^+_{2k}$$ to be the subset of
classes $z$ which can be represented by a positive holomorphic $k$-cycle.

\Theorem{\MH.1}  {\sl Suppose Conjecture B holds. Then the cycle $\G$ bounds
a positive holomorphic $p$-chain  of mass $\leq \Lambda$ in $X$  if and only
if any of the following conditions holds:
\medskip

(a)\ \ 
$\widetilde{\Link}_{\lambda}(\G,Z) \geq -\Lambda$ \ \ for all  
positive algebraic cycles  $Z\subset X-\G$  of codimension-$p$
 
\qquad with cohomology
class  $\ell[\omega]^p$ for $\ell\in \bbz^+$.
\medskip

(b) \ \ 
$\widetilde{\Wind}_{\lambda}(\G,\a) \geq -\Lambda$ \ \ for all  
sparks $\a$  satisfying equation (\MH.1) with $Z$ as above.

\medskip

(c) \ \ 
${1\over  p!} \int_\G\b \geq -\Lambda$ \ \ for all smooth forms $\b\in\ce^{2p-1}(X)$ for which
$d^{p,p}\b +\o^p \geq 0$ is 

\qquad weakly positive on $X$.

\medskip

(d) \ \ 
There exists $\tau\in H_{2p}(X, |\G|;\bbz)$ with  $\partial \tau = [\G]$ such that 
$
\tau\bullet [Z] \ \geq\ 0$ for all 

\qquad \  $[Z]$ as above and $\Lambda = (\tau, {1\over p!} [\omega^p])$
}

\pf  Condition (a) represents Theorem \PM.2 above.  Conditions (a) and (b) are obviously equivalent by equation (\MH.2).
To check  the necessity of Condition (c) suppose that $\G=dT$ where $T$ is
 a positive holomorphic $p$-chain. Then one has $\int_\G \b = \int_T d\b  = \int_T d^{p,p}\b
 = \int_T (d^{p,p}\b +\o^p-\o^p) \geq -\int_T \o^p = -p! \M(T)$. On the other hand, Condition
 (c) implies Condition (a) since the subvarieties $Z$ in question arise by intersection with
 subvarieties $\widetilde Z$ in the ambient projective space, where we can mollify to obtain
 smooth $(p-1,p-1)$-forms $\a_\epsilon$ with $d\a_\epsilon = \o^p-{\widetilde Z}_\epsilon$
and ${\widetilde Z}_\epsilon \to {\widetilde Z}$.

For Condition (d), suppose there exists $T$ as above.  
Then the class $\tau = [T] \in H_{2p}(X, |\G|;\bbz)$  has the stated properties.
Conversely, given $\tau$, choose any 2$p$-chain $N\in \tau$. Then $dN=\G$ and
for any $Z$ as above we have $0\leq \tau\bullet [Z] = N\bullet Z =
N\bullet Z -\ell\int_N\o^p+\ell\int_N\o^p = \ell p!\, \widetilde{\Link}_{\lambda}(\G,Z) + \ell \tau(\o^p)$.
\qed

\medskip
\Cor{\MH.2}  {\sl  If
$H^+_{2(n-p), {\rm alg}}$ is contained in a proper subcone of
$H^+_{2(n-p)}$, that is, if there exists $\tau\in H_{2p}(X,|\G|;\bbz)$ with 
$d\tau\neq 0$ and $\tau\bullet u\geq0$ for $u\in H^+_{2(n-p), {\rm alg}}$,
then there exists a positive holomorphic $p$-chain $T$ on $X$ with 
non-empty boundary supported in $|\G|$.}
\medskip


This question of holomorphic representability is discussed in detail in
[HL$_5$].

\vskip .3in



\centerline{\bf References}

\vskip .2in

\noindent
[AW]  H. Alexander and J. Wermer, {\sl Linking numbers
and boundaries of varieties}, Ann. of Math.
{\bf 151} (2000),   125-150.

 \smallskip

\noindent
[DH]   P. Dolbeault and G. Henkin,   
{\sl Cha\^ines holomorphes de bord donn\'e dans $\bbc\bbp^n$},    
Bull.  Soc. Math. de France,
{\bf  125}  (1997), 383-445.

\smallskip

\noindent
[F]   H. Federer, Geometric Measure  Theory,
 Springer--Verlag, New York, 1969.

 \smallskip

\noindent
[Ha]   B. Harris, {\sl
Differential characters and the Abel-Jacobi map},
pp. 69-86 in ``Algebraic K-theory;  Connections with Geometry and
Topology'', Jardine and Snaith (eds.), Kluwer Academic Publishers,
1989.
\smallskip

\noindent
[H]  F.R. Harvey,
Holomorphic chains and their boundaries, pp. 309-382 in ``Several Complex
Variables, Proc. of Symposia in Pure Mathematics XXX Part 1'', 
A.M.S., Providence, RI, 1977.

 \noindent 
\noindent
[HL$_1$] F. R. Harvey and H. B. Lawson, Jr, {\sl On boundaries of complex
analytic varieties, I}, Annals of Mathematics {\bf 102} (1975),  223-290.

 \smallskip

 \noindent
[HL$_2$] F. R. Harvey and H. B. Lawson, Jr, {\sl Boundaries of 
varieties in projective manifolds}, J. Geom. Analysis,  {\bf 14}
no. 4 (2005), 673-695. ArXiv:math.CV/0512490.
 \smallskip

 \noindent 
\noindent
[HL$_3$] F. R. Harvey and H. B. Lawson, Jr, {\sl Projective hulls and
the projective Gelfand transformation}, Asian J. Math. {\bf 10}, no. 2 (2006), 279-318. ArXiv:math.CV/0510286.

 \smallskip

\noindent
[HL$_4$] F. R. Harvey and H. B. Lawson, Jr, {\sl Projective linking and boundaries of positive holomorphic chains in projective manifolds, Part I},   ArXiv:math.CV/0512379

 \smallskip

\noindent
[HL$_5$] F. R. Harvey and H. B. Lawson, Jr, {\sl Relative holomorphic cycles and duality}, Preprint, Stony Brook, 2006.

 \smallskip

\noindent
[HLW] F. R. Harvey, H. B. Lawson, Jr. and J. Wermer, {\sl On the projective hull of certain curves in $\bbc^2$},  Stony Brook Preprint, 2006.

 \smallskip

\noindent
[HLZ] F. R. Harvey, H. B. Lawson, Jr. and J. Zweck, {\sl A
deRham-Federer theory of differential characters and character duality},
Amer. J. of Math.  {\bf 125} (2003), 791-847.

 \smallskip

   \noindent
[HS]    F.R. Harvey and B. Shiffman,    {\sl  A characterization of
holomorphic chains},    Ann. of Math.,
 {\bf 99}  (1974), 553-587.

\smallskip

\end

\Theorem{\PHC.3} {\sl Let $\lambda$ be a holomorphic line
bundle with positive curvature form $\omega$ on a compact complex
manifold $X$. Then functions of the form $\log\|\s\|^{1\over{\ell}}$ for
$\s\in  H^0(X, \co{}(\lambda^{\ell}))$  are weakly dense in
$\psh(X)$.  }

\vskip .3in\
\noindent
{\bf First argument:}   Let $\chi_\e$ be the characteristic funtion of the
complement of an $\e$-tubular neighborhood of  $\G$, and set $T_\e = \chi_\e
T$. Then  $\M(T_\e) \leq \L$  for all $\e$, and by Federer slicing theory [F],
we have $\M(d(T_\e) = \M(d\chi_\e\w T) < \infty$ for almost all $\e$, and in
particular for some sequence $\e_j\to 0$.

Consider $T_0\equiv \lim_{\e\to 0} T_\e  = \lim_{j\to \infty} T_{\e_j}$.
Since $T_{\e_j}$ is normal and $F(T_{\e_j}-T_0)\leq  \M(T_{\e_j}-T_0)  
\arr 0$ as $j\to\infty$, where $F$ denotes the flat norm, we conclude that
$T_0$ is flat. Hence $dT_0$ is flat, and so 
$$
d(T+R-T_0)\ =\ \G-dT_0  \qquad {\rm is \ flat.}
\eqno(\PAW.6)
$$
Now by (\PAW.5) we have $\supp (T-T_0+R)\subset \G$. Let $\pi$ be a smooth
retraction of a tubular neighborhood of $\G$ onto $\G$.   Then
$$
\G-dT_0 \ =\ \pi_* (\G-dT_0) \ =\ d  \pi_*(T+R-T_0)\ = \ 0
$$
where the first equality comes from (\PAW.6) and the Federer Flat
Support Theorem [F, page 372], and the last equality holds for formal degree
reasons.

\bigskip
{\bf REESE, We still need to prove that the $c_j$'s are integers.}